\input amstex
\documentstyle{amsppt}
\topmatter
\magnification=\magstep1
\pagewidth{5.2 in}
\pageheight{6.7 in}
\abovedisplayskip=10pt
\belowdisplayskip=10pt
\NoBlackBoxes
\title
Note on  Dedekind type DC sums
\endtitle
\author  Taekyun Kim  \endauthor
\affil\rm{{Division of General Education-Mathematics,}\\
{ Kwangwoon University, Seoul 139-701, S. Korea}\\
{e-mail: tkkim$\@$kw.ac.kr}\\\\
}\endaffil

\abstract{In this paper we study the Euler polynomials and functions
and derive some interesting formulae related to the Euler
polynomials and functions. From those formulae we consider Dedekind
type DC(Daehee-Changhee)sums  and prove reciprocity laws related to
DC sums.}
\endabstract
\thanks 2000 Mathematics Subject Classification  11S80, 11B68 \endthanks
\thanks Key words and phrases: Dedekind sum, Hardy sum, Euler polynomials, DC sum \endthanks
\rightheadtext{   } \leftheadtext{Note on Dedekind type DC sums  }
\endtopmatter

\document

\head 1. Introduction/Preliminaries \endhead

The Euler numbers are defined as
$$ \frac{2}  {e^t +1}=\sum_{n=0}^{\infty}E_n \frac{t^n}{n!}, \text{ $|t|< \pi$}, \text{ (see [1-31])},\tag1$$
and the Euler polynomials are also defined as
$$\frac{2}{e^t +1}e^{xt}=\sum_{n=0}^{\infty}E_n(x)\frac{t^n}{n!},
\text{  $|t|<\pi$, ( see [4, 5, 6] ).} \tag2$$ The first few of
Euler numbers  are $1, -\frac{1}{2}, 0, \frac{1}{4}, $ and
$E_{2k}=0$ for $k=1, 2, 3, \cdots.$ From (1) and (2), we can
easily derive the following.
$$ E_n(x)=\sum_{l=0}^n \binom{n}{l}E_lx^{n-l}, \text{ where
$\binom{n}{l}=\frac{n(n-1)\cdots (n-l+1)}{l!}$, (see [4, 5, 6]).}
\tag3$$ In (1), it is easy to see that $E_0=1$,
$E_n(1)+E_n=2\delta_{0,n},$
 where $\delta_{0,n}$ is the Kronecker symbol. That is,
 $E_n(1)=-E_n$ for $n=1, 2, 3, \cdots,$ (see [4, 5]). We denote $\bar
 E_n(x)$ the $n$-th Euler function given by  the Fourier expansion.
 $$\bar E_n(x)=m!2\sum_{n=-\infty}^{\infty} \frac{e^{(2n+1)\pi
 ix}}{((2n+1)\pi i)^{m+1}}, \text{ (see [4,5])}, $$
 which, for $0\leq x<1,$ reduces the $n$-th Euler polynomials.

 By (3), we easily see that
 $$\frac{dE_n(x)}{dx}=\frac{d}{dx}\sum_{k=0}^{n}\binom{n}{k}=n\sum_{k=0}^n\binom{n-1}{k}E_kx^{n-1-k}=nE_{n-1}(x).\tag4$$
From (4), we note that
$$\int_{0}^{x}E_n(t)dt=\frac{1}{n+1}E_{n+1}(x), \text{ (see [4])}. \tag5$$
By the definitions of the Euler numbers and the Euler polynomials,
we easily see that
$$2\sum_{k=0}^{n-1}(-1)^ke^{kt}=2\frac{(-1)^ne^{nt}+1}{e^t
+1}=\sum_{l=0}^{\infty}\left((-1)^nE_l(n)+E_l\right)\frac{t^l}{l!}.\tag6$$
Thus, we have
$$2\sum_{k=0}^{n-1}(-1)^kk^l=(-1)^nE_l(n)+E_l.\tag7$$

It is well known that the classical Dedekind sums $S(h,k)$ first
arose in the transformation formula of the logarithm of Dedekind
eta-function (see [17, 25, 26, 28]). If $h$ and $k$ are relative
prime integers with $k>0$, then  Dedekind  sum is defined as
$$S(h, k)=\sum_{u=1}^{k-1}(( \frac{u}{k}))((\frac{hu}{k})), \text{ (see [17, 25, 26, 28])}\tag8$$
where $((x))$ is defined as
$$\aligned
((x))&= x-[x]-\frac{1}{2}, \text{ if $x$ is not an integer, }\\
      &=0, \text{ otherwise, }
 \endaligned$$
where $[x ]$ is the largest integer $\leq x$, (cf. [17]).

Generalized Dedekind sums $S_p(h, k)$ are defined as
$$ S_p(h, k)=\sum_{a=1}^{k-1}\frac{a}{k}\bar B_p(\frac{ah}{k}),
\text{ (see [1, 17, 18, 25, 26, 28, 31]),}\tag9$$ where $h$ and
$k$ are relative prime positive integers and $\bar B_p(x)$ are the
$p$-th Bernoulli functions, which are defined as
$$\bar B_p(x)=B_p(x-[x])=-p!(2 \pi i)^{-p}\sum_{m=-\infty, m\neq
0}^{\infty}m^{-p}e^{2\pi i nx}, \text{ (see [17, 25, 26])},$$
where $B_p(x)$ are the $p$-th ordinary Bernoulli polynomials.

Recently Y. Simsek have studies $q$-Dedekind type sums related to
$q$-zeta function and basic $L$-series (see [31, 18]). He also
studies $q$-Hardy-Berndt type sums associated with $q$-Genocchi
type zeta and $q$-$l$-functions related to previous author's paper
(see [18, 31]). In this paper we consider Dedekind type
DC(Daehee-Changhee) sums as follows.
$$T_p(h, k)=2\sum_{u=1}^{k-1}(-1)^{u-1}\frac{u}{k}\bar
E_p(\frac{hu}{k} ), \text{ ( $h\in \Bbb Z_+$ )}, $$ where $\bar
E_p(x)$ are the $p$-the Euler functions.  Note that $T_p(h, k)$ is
the similar form of generalized Dedekind type sums. Finally, we
prove the following reciprocity law for an odd $p$:
$$\aligned
&k^pT_p(h, k)+h^pT_p(k,h)\\
&=2\sum_{\Sb u=0\\ u-[\frac{hu}{k}]\equiv 1\mod
2\endSb}^{k-1}\left(kh(
E+\frac{u}{k})+k(E+h-[\frac{hu}{k}])\right)^p
+\left(hE+kE\right)^p+(p+2)E_p,
\endaligned$$
where $h, k$ are relative prime  positive integers and
$$(Eh+Ek)^{n+1}=\sum_{l=0}^{n+1}\binom{n+1}{l}E_lh^lE_{n+1-l}k^{n+1-l}.$$

\head 2. On the reciprocity law for Dedekind type DC sums
\endhead

In this section, we assume $p \in \Bbb N$ with  $p\equiv 1 \mod 2$.
By the definition of the Euler polynomials, we see that

$$\aligned
E_p(x+y)&=\sum_{s=0}^p\binom{p}{s}(x+y)^{p-s}E_s=\sum_{s=0}^p\binom{p}{s}E_s\sum_{k=0}^{p-s}\binom{p-s}{k}
x^ky^{p-s-k}\\
&=\sum_{s=0}^p\binom{p}{s}\sum_{j=0}^s\binom{s}{j}E_jx^{s-j}y^{p-s}=\sum_{s=0}^p\binom{p}{s}E_s(x)
y^{p-s}.
\endaligned\tag10$$
From (2), we can also derive
$$E_p(mx)=m^p\sum_{s=0}^{m-1}E_p(x+\frac{s}{m})(-1)^s. \tag11$$
By (5), we easily see that
$$\int_0^1 x E_p(x)dx
=\frac{E_{p+1}(1)}{p+1}-\frac{E_{p+1}(1)}{p+1}+\frac{E_{p+1}}{p+1}=\frac{E_{p+1}}{p+1}=0,\tag12$$
and
$$\int_0^1 x E_p(x) dx=\sum_{s=0}^p \binom{p}{s}E_s\int_0^1
x^{p-s+1}dx=\sum_{s=0}^{p}\binom{p}{s}\frac{E_s}{p-s+2}. \tag13$$ By
(12) and (13), we obtain the following lemma.

\proclaim{ Lemma 1} For $p\in\Bbb N$ with $p\equiv 1 \mod 2$, we
have
$$\sum_{s=0}^p
\binom{p}{s}\frac{E_s}{p-s+2}=\frac{E_{p+1}}{p+1}=0.$$
\endproclaim
For $s\in\Bbb N$ with $s\equiv 0 \mod 2$ and $s<p$, we have
$$\frac{d^s(xE_p(x))}{(dx)^s}\big|_{x=1}=s!\binom{p}{s}E_{p-s}
(1)=-s!\binom{p}{s}E_{p-s} ,\tag14$$

and, from (3), we note that
$$\frac{d^s(xE_p(x))}{(dx)^s}\big|_{x=1}=s!\sum_{v=0}^{p-s}\binom{p-v+1}{s}\binom{p}{v}E_v
=s!\sum_{v=0}^{p}\binom{p}{v}\binom{p-v+1}{s}E_v.\tag15$$

 By (14) and (15), we obtain the following theorem.

 \proclaim{Theorem 2} For $s \in \Bbb N$ with $s\equiv 0 \mod 2$ and
 $ s>p$, we have
 $$\sum_{v=0}^p
 \binom{p}{v}\binom{p-v+1}{s}E_v=-\binom{p}{s}E_{p-s}=\binom{p}{s}E_{p-s}(1).\tag16$$
 \endproclaim
Let us define  Dedekind type DC sums as follows.
$$T_p(h, k)=2\sum_{u=1}^{k-1}(-1)^{u-1}\frac{u}{k}\bar
E_p(\frac{hu}{k}), \text{ ( $h\in\Bbb Z_{+}$),}\tag17$$ where $\bar
E_p(x)$ is the $p$-th Euler function.

For $m\equiv 1 \mod 2$, we have
$$\aligned
T_p(1, m)&=2\sum_{u=1}^{m-1}(-1)^{u-1}\frac{u}{m}\sum_{v=0}^p
\binom{p}{v}E_v\left(\frac{u}{m}\right)^{p-v}\\
&=\sum_{v=0}^p\binom{p}{v}E_v
m^{-(p+1-v)}2\sum_{u=1}^{m-1}(-1)^{u-1}u^{p-v+1}.
\endaligned\tag18$$

By (7) and (18), we obtain the following theorem.

\proclaim{ Theorem 3} For $ m\equiv 1 \mod 2,$ we have
$$T_{p}(1, m)=\sum_{v=0}^p \binom{p}{v}E_v m^{-(p+1-v)}
\left(E_{p-v+1}(m)-E_{p-v+1} \right). \tag19$$
\endproclaim

From (3) we can also derive
$$E_{p-v+1}(m)-E_{p-v+1}=\sum_{i=0}^{p-v}\binom{p-v+1}{i} m^{p+1-v-i}E_i \tag20$$
so that we find
$$\aligned
T_p(1, m)&=\sum_{v=0}^p\binom{p}{v}m^{-(p+1-v)}E_v\sum_{i=0}^{p-v}\binom{p-v+1}{i}m^{p+1-v-i}E_i\\
&=\frac{1}{m^p}\sum_{v=0}^p\binom{p}{v}E_v\sum_{i=0}^{p-v}\binom{p-v+1}{i}E_i
m^{p-i}.
\endaligned\tag21$$
Therefore, we obtain the following corollary.

\proclaim{ Corollary 4} For $m\equiv 1 \mod 2$, we have
$$m^p T_p(1, m)=\sum_{v=0}^p
\binom{p}{v}E_v\sum_{i=0}^{p-v}\binom{p-v+1}{i}E_im^{p-i}. \tag22$$
\endproclaim

Interchanging the order of summation in (22), we obtain
$$\aligned
m^pT_p(1,
m)&=\sum_{i=0}^p\sum_{v=0}^{p-i}\binom{p}{v}E_v\binom{p-v+1}{i}E_i
m^{p-i}=\sum_{i=1}^{p-2}\sum_{v=0}^{p-i}\binom{p}{v}E_v\binom{p-v+1}{i}E_im^{p-i}\\
&+\binom{p+1}{p}E_p+\sum_{v=0}^p\binom{p}{v}E_vm^p
+\sum_{v=0}^1\binom{p}{v}E_v\binom{p-v+1}{p-1}E_{p-1}m\\
&=\sum_{i=1}^{p-2}\sum_{v=0}^{p-i}\binom{p}{v}E_v\binom{p-v+1}{i}E_im^{p-i}
+(p+1)E_p+\sum_{v=0}^p\binom{p}{v}E_vm^p.
\endaligned$$
Therefore, we obtain the following proposition.

\proclaim{ Proposition 5} For $m\in\Bbb N$ with $m\equiv 1 \mod
2$, we have
$$m^pT_p(1, m)=\sum_{v=0}^p\binom{p}{v}E_vm^p+\sum_{i=1}^{p-2}\sum_{v=0}^{p-i}\binom{p}{v}E_v\binom{p-v+1}{i}E_im^{p-i}
+(p+1)E_p .\tag23$$
\endproclaim
In the sum over $i$, the only non-vanishing terms are those for
which the index $i$ is odd. Hence, since $i<p$ in this sum we may
use (3) and Theorem 2 to obtain
$$\aligned
m^pT_p(1, m)&=m^pE_p(1)+(p+1)E_p
+\sum_{i=1}^{p-2}\binom{p}{i}E_{p-i}(1)E_im^{p-i}\\
&=\sum_{i=0}^p\binom{p}{i}E_{p-i}(1)E_im^{p-i}+pE_p .
\endaligned\tag24$$
Therefore, we obtain the following theorem.

\proclaim{Theorem 6} For odd $p$ with $p>1$, $m\in\Bbb Z_{+}$ with
$ m\equiv 1 \mod 2$, we have

$$m^pT_p(1, m)=\sum_{i=0}^p \binom{p}{i}E_{p-i}(1)E_i
m^{p-i}+pE_p.$$
\endproclaim
 Now we employ the symbolic notation as $E_n(x)=(E+x)^n$.
It is easy to show that
$$\aligned
&k^p\sum_{u=0}^{k-1}(-1)^u\sum_{s=0}^p\binom{p}{s}h^sE_s(\frac{u}{k})E_{p-s}(h-[\frac{hu}{k}])
=k^p\sum_{u=0}^{k-1}(-1)^u\left(h(E+\frac{u}{k})+(E+h-[\frac{hu}{k}])
\right)^p\\
&=k^p\sum_{u=0}^{k-1}(-1)^u\left(Eh+E+h+\frac{1}{2}-\frac{1}{2}+huk^{-1}-[\frac{hu}{k}]\right)^p\\
&=k^p\sum_{u=0}^{k-1}(-1)^u\left(Eh+E+h+\frac{1}{2}+\bar
E_1(\frac{hu}{k})\right)^p .
\endaligned\tag25$$

Now as the index $u$ range through the values $u=0, 1, 2, \cdots,
k-1,$ the product $hu$ range through a complete residue system
modulo $k$ since $(h, k)=1$ and due to the periodically of $\bar
E_1(x)$, the term $\bar E_1(\frac{hu}{k})$ may be replaced $\bar
E_1(\frac{u}{k})$ without alternating  the sum over $u$. For
$k\in\Bbb Z_{+}$ with $k\equiv 1 \mod 2$, we have
$$\aligned
&(25)=k^p\sum_{u=0}^{k-1}(-1)^u\left(E+Eh+h+\frac{1}{2}+\bar
E_1(\frac{u}{k})\right)^p
=k^p\sum_{u=0}^{k-1}(-1)^u\left((E+\frac{u}{k})+h(E+1)\right)^p\\
&=k^p\sum_{u=0}^{k-1}(-1)^u\sum_{s=0}^p\binom{p}{s}E_s(\frac{u}{k})h^{p-s}E_{p-s}(1)\\
&=\sum_{s=0}^p\binom{p}{s}k^{p-s}\left(k^s\sum_{u=0}^{k-1}(-1)^uE_s(\frac{u}{k})\right)h^{p-s}E_{p-s}(1)
=\sum_{s=0}^p\binom{p}{s}k^{p-s}E_sh^{p-s}E_{p-s}(1).
\endaligned$$
Therefore, we obtain the following theorem.

\proclaim{ Theorem 7} Let $h, k$ be natural numbers with $(h,
k)=1$. For odd $p$ with $p>1$, and $k\equiv 1 \mod 2$, we have
$$\sum_{s=0}^p\binom{p}{s}k^{p-s}E_sh^{p-s}E_{p-s}(1)
=k^p\sum_{u=0}^{k-1}(-1)^u\sum_{s=0}^p\binom{p}{s}h^sE_s(\frac{u}{k})
E_{p-s}(h-[\frac{hu}{k}]).$$
\endproclaim
Let $T$ be the sum of
$$\aligned
T&=k^pT_p (h, k)+h^pT(k,h)\\
&=2k^p\sum_{u=1}^{k-1}(-1)^{u-1}\frac{u}{k}\bar E_p(\frac{hu}{k})
+2h^p\sum_{v=0}^{h-1}(-1)^{v-1}\frac{v}{h}\bar E_p(\frac{kv}{h}).
\endaligned\tag26$$
We assume first that $p>1$ and $h, k \in \Bbb N$ with $h\equiv 1
\mod 2$, and $k\equiv 1 \mod 2$.
$$\bar E_p(\frac{h}{k}u)=h^p\sum_{v=0}^{h-1}(-1)^v\bar
E_p(\frac{u}{k}+\frac{v}{h}), \tag27$$ and
$$\bar E_p(k\frac{v}{h})=k^p\sum_{u=0}^{k-1}(-1)^u\bar
E_p(\frac{v}{h}+\frac{u}{k}).$$

From (26) and (27), we can easily derive the following (28).

$$\aligned
T&=(hk)^p2\sum_{u=1}^{k-1}(-1)^{u-1}\frac{u}{k}\sum_{v=0}^{h-1}(-1)^v
\bar E_p(\frac{u}{k}+\frac{v}{h})\\
&+(hk)^p2\sum_{v=1}^{h-1}(-1)^{v-1}\frac{v}{h}\sum_{u=0}^{k-1}(-1)^u
\bar E_p(\frac{v}{h}+\frac{u}{k})\\
&=(hk)^p2\sum_{u=0}^{k-1}\sum_{v=0}^{h-1}(-1)^{u+v-1}\left(\frac{uh+vk}{hk}\right)E_p(\frac{u}{k}+\frac{v}{h}).
\endaligned\tag28 $$
Therefore, we obtain the following theorem.

\proclaim{Theorem 8} Let $h, k \in \Bbb N$ with $h\equiv 1 \mod 2$
and $k\equiv 1 \mod 2$. For $p>1$, we have
$$k^pT_p(h, k)+h^pT(k,
h)=2(hk)^p\sum_{u=0}^{k-1}\sum_{v=0}^{h-1}(-1)^{u+v-1}(uh+vk)(hk)^{-1}E_p(\frac{u}{k}+\frac{v}{h}).$$
\endproclaim
Now as the indices $u$ and $v$ run through the range $u=0, 1, 2,
cdots, k-1$, $v=0, 1, 2, \cdots, h-1 ,$ respectively, the linear
combination $uh+vk$ ranges through a complete residue system
modulo $hk ,$ and each term $uh +vk$ satisfies the inequalities
$0\leq uh +vk <2hk.$ If we define the sets
$$\aligned
& A=\{uh+vk| 0\leq uh+vk< hk \},  B=\{uh+vk|hk+1 \leq uh+vk
<2hk\},\\
& C=\{ \lambda |0 \leq \lambda \leq hk-1 \} .
\endaligned$$
Let $h, k\in \Bbb N$ with $h\equiv 1 \mod 2$ and $k\equiv 1 \mod
2$. From (28), we note that

$$T=(hk)^p\left(2\sum_{\lambda \in
A}\frac{\lambda}{hk}(-1)^{\lambda-1}\bar E_p(\frac{\lambda}{hk})
+2\sum_{\lambda\in B}\frac{\lambda}{hk}(-1)^{\lambda-1}\bar
E_p(\frac{\lambda}{hk})\right). \tag29$$

Now if $y\in B$, then $y=hk+\lambda$, where $\lambda \in \Bbb C$,
but $\lambda \notin A$ (for if $\lambda \in A$ then we have
$\lambda \equiv y \mod hk$), but $A \bigcup B$ forms a complete
residue system modulo $hk$. Hence, we have

$$2\sum_{y\in B}\frac{y}{hk}(-1)^{y-1}\bar E_p(\frac{y}{hk})=2
\sum_{\lambda \in C\setminus A} (-1)^{\lambda -1}\bar
E_p(\frac{\lambda}{hk})+2\sum_{\lambda\in C\setminus A}
\frac{\lambda}{hk}(-1)^{\lambda-1}\bar E_p(\frac{\lambda}{hk}).
\tag30$$

By (29) and (30), we see that

$$\aligned
T&=(hk)^p\big\{2\sum_{\lambda \in
A}\frac{\lambda}{hk}(-1)^{\lambda-1}\bar E_p
(\frac{\lambda}{hk})+2\sum_{\lambda\in C\setminus
A}(-1)^{\lambda-1}\bar E_p(\frac{\lambda}{hk})
\\
&+2\sum_{\lambda\in C\setminus A}
\frac{\lambda}{hk}(-1)^{\lambda-1}\bar
E_p(\frac{\lambda}{hk})\big\}\\
&=(hk)^p\big\{2\sum_{\lambda=0}^{hk-1}(-1)^{\lambda-1}\frac{\lambda}{hk}\bar
E_p(\frac{\lambda}{hk})+2\sum_{\lambda=0}^{hk-1}(-1)^{\lambda-1}\bar
E_p(\frac{\lambda}{hk})\\
&-2\sum_{\Sb u=0\\ 0\leq uh+vk<hk
\endSb}^{k-1}\sum_{v=0}^{h-1}(-1)^{u+v-1}\bar
E_p(\frac{uh+vk}{hk})\big\}.
\endaligned$$
It is easy to see that
$$2\sum_{\lambda =0}^{hk-1}\bar
E_p(\frac{\lambda}{hk})(-1)^{\lambda-1}=2(hk)^{-p}\bar
E_p(0)=2(hk)^{-p}E_p.$$ Hence, we have
$$T=(hk)^p\left(T_p(1, kh)+2(hk)^{-p}E_p-S\right), \tag31$$
where
$$\aligned
S&=2\sum_{\Sb 0\leq u \leq k-1\\ 0\leq uh+vk<hk \endSb}\sum_{0\leq
v \leq h-1} (-1)^{u+v-1} \bar E_p(\frac{uh+vk}{hk})\\
&=2\sum_{\Sb 0\leq u\leq k-1\\ 0\leq
\frac{u}{k}+\frac{v}{h}<1\endSb}\sum_{0\leq v \leq
h-1}(-1)^{u+v-1}E_p(\frac{u}{k}+\frac{v}{h}).\endaligned$$
 From the definition of $S$, we note that
 $$\aligned
&S=2\sum_{u=0}^{k-1}\sum_{v=0}^{[h-\frac{hu}{k}]}(-1)^{u+v-1}E_p(\frac{u}{k}+\frac{v}{h})
=2\sum_{u=0}^{k-1}\sum_{v=0}^{[h-\frac{hu}{k}]}(-1)^{u+v-1}\left(E+\frac{u}{k}+\frac{v}{h}\right)^p\\
&=2\sum_{s=0}^p\binom{p}{s}h^{s-p}\sum_{u=0}^{k-1}(-1)^{u-1}E_s(\frac{u}{k})\sum_{v=0}^{[h-\frac{hu}{k}]}(-1)^vv^{p-s}\\
&=\sum_{s=0}^p
\binom{p}{s}h^{s-p}\sum_{u=0}^{k-1}(-1)^{u-1}E_s(\frac{u}{k})\left(2\sum_{v=0}^{h-1-[\frac{hu}{k}]}(-1)^vv^{p-s}\right)\\
&=\sum_{s=0}^p\binom{p}{s}h^{s-p}\sum_{u=0}^{k-1}(-1)^{u-1}E_s(\frac{u}{k})
\left((-1)^{h-[\frac{hu}{k}]}E_{p-s}(h-[\frac{hu}{k}])+E_{p-s}\right)\\
&=\sum_{s=0}^p\binom{p}{s}h^{s-p}\sum_{u=0}^{k-1}(-1)^{u-[\frac{hu}{k}]}E_s(\frac{u}{k})E_{p-s}(h-[\frac{hu}{k}])
+\sum_{s=0}^p\binom{p}{s}h^{s-p}\sum_{u=0}^{k-1}(-1)^{u-1}E_s(\frac{u}{k})E_{p-s}\\
&=\sum_{s=0}^p\binom{p}{s}h^{s-p}\sum_{u=0}^{k-1}(-1)^{u-[\frac{hu}{k}]}E_s(\frac{u}{k})E_{p-s}(h-[\frac{hu}{k}])
-h^{-p}\sum_{s=0}^p\binom{p}{s}h^sk^{-s}E_{p-s}E_s.
 \endaligned$$
Returning to (31), we have
$$\aligned
T&=(hk)^p\big\{T_p(1,
kh)+2(kh)^{-p}E_p-\sum_{s=0}^{p}\binom{p}{s}h^{s-p}\sum_{u=0}^{k-1}(-1)^{u-[\frac{hu}{k}]}
E_{p-s}(h-[\frac{hu}{k}])E_s(\frac{u}{k})\\
&+h^{-p}\sum_{s=0}^p \binom{p}{s}h^sk^{-s}E_{p-s}E_s\big\}.
\endaligned$$
By Theorem 6 we see that
$$\aligned
T&=\sum_{s=0}^p\binom{p}{s}E_sE_{p-s}(1)(hk)^{p-s}-\sum_{s=0}^p\binom{p}{s}h^sk^p\sum_{u=0}^{k-1}(-1)^{u-[\frac{hu}{k}]}
E_s(\frac{u}{k})E_{p-s}(h-[\frac{hu}{k}])\\
&+(p+2)E_p +\sum_{s=0}^p\binom{p}{s}h^sk^{p-s}E_s E_{p-s}.
\endaligned$$
From Theorem 7, we can also derive  the following equation (32).

$$\aligned
T&=\sum_{s=0}^p\binom{p}{s}k^ph^s\sum_{u=0}^{k-1}E_s(\frac{u}{k})E_{p-s}(h-[\frac{hu}{k}])(1-(-1)^{u-[\frac{hu}{k}]})\\
&+\sum_{s=0}^p \binom{p}{s}h^sk^{p-s}E_sE_{p-s}+(p+2)E_p
 \endaligned\tag32$$
Therefore, we obtain the following theorem.

 \proclaim{ Theorem 9}
Let $h, k \in \Bbb N$ with $h\equiv 1 \mod 2$ and $k\equiv 1 \mod
2$ and let $(h, k)=1$. For $p>1$, we have
$$\aligned
&k^pT_p(h, k)+h^pT_p(k, h) \\
&=2\sum_{\Sb u=0\\u-[\frac{hu}{k}]\equiv 1 (\mod 2)
\endSb}^{k-1}\left(kh(E+\frac{u}{k})+k(E+h-[\frac{hu}{k}])\right)^p
+(hE+kE)^p+(p+2)E_p,\endaligned$$ where
$$ (hE+kE)^p=\sum_{s=0}^p\binom{p}{s}h^sE_sk^{p-s}E_{p-s}.$$
\endproclaim

\enddocument